\documentclass[12pt]{amsart}
\usepackage{amscd,amsmath,amssymb,amsfonts}
\theoremstyle{plain}
\newtheorem{thm}{Theorem}
\newtheorem{lem}[thm]{Lemma}

\newtheorem{prop}[thm]{Proposition}
\newtheorem{remark}[thm]{Remark}

\newtheorem{defn}[thm]{Definition}

\newtheorem{conj}[thm]{Conjecture}
\numberwithin{thm}{section}
\numberwithin{equation}{section}

\newcommand{\al}{\alpha}

\newcommand{\eq}[2]{\begin{equation}\label{#1}#2 \end{equation}}
\newcommand{\ml}[2]{\begin{multline}\label{#1}#2 \end{multline}}

\newcommand{\inj}{\hookrightarrow}

\newcommand{\Ker}{{\rm Ker}}
\newcommand{\Pic}{{\rm Pic}}

\newcommand{\Spec}{{\rm Spec \,}}

\newcommand{\Tr}{{\rm Tr}}

\newcommand{\rt}{{\rm res \ Tr}}

\newcommand{\be}{{\beta}}
\newcommand{\ta}{{\theta}}
\newcommand{\ga}{{\gamma}}

\newcommand{\sD}{{\mathcal D}}
\newcommand{\sE}{{\mathcal E}}

\newcommand{\sG}{{\mathcal G}}

\newcommand{\sK}{{\mathcal K}}
\newcommand{\sL}{{\mathcal L}}

\newcommand{\sO}{{\mathcal O}}
\newcommand{\sP}{{\mathcal P}}

\newcommand{\sU}{{\mathcal U}}

\newcommand{\A}{{\mathbb A}}

\newcommand{\C}{{\mathbb C}}

\newcommand{\G}{{\mathbb G}}
\renewcommand{\H}{{\mathbb H}}

\newcommand{\N}{{\mathbb N}}
\renewcommand{\P}{{\mathbb P}}
\newcommand{\Q}{{\mathbb Q}}
\newcommand{\R}{{\mathbb R}}

\newcommand{\Z}{{\mathbb Z}}

\begin{document}

\title[Irregular Connections]{Gau\ss-Manin Determinant Connections and Periods for
Irregular Connections} 
\author{Spencer Bloch}
\address{Dept. of Mathematics,
University of Chicago,
Chicago, IL 60637,
USA}
\email{bloch@math.uchicago.edu}

\author{H\'el\`ene Esnault}
\address{Mathematik,
Universit\"at Essen, FB6, Mathematik, 45117 Essen, Germany}
\email{esnault@uni-essen.de}
\date{December 13, 1999}
\begin{abstract} Gau\ss-Manin determinant connections associated to irregular
connections on a curve are studied. The determinant of the
Fourier transform of an 
irregular connection is calculated. The determinant
of cohomology of the standard rank 2 Kloosterman sheaf is computed
modulo 2 torsion.
Periods associated to
irregular connections are 
studied in the very basic $\exp(f)$ case, and analogies with the Gau\ss-Manin
determinant are discussed. 
\end{abstract}
\maketitle
\begin{quote}Everything's so awful reg'lar a body can't stand it.\\
\\
The Adventures of Tom Sawyer \\
Mark Twain
\end{quote}

\section{Introduction}

A very classical area of mathematics, at the borderline between applied
mathematics, algebraic geometry, analysis, mathematical physics and number
theory is the theory of systems of linear differential equations
(connections). There is a vast literature focusing on regular singular
points, Picard-Fuchs differential
equations, Deligne's Riemann-Hilbert correspondence and its extension to
$\sD$-modules, and more recently various index theorems in geometry.

On the
other hand, there are some very modern themes in the subject which remain
virtually untouched. On the arithmetic side, one may ask, for example, how
irregular connections can be incorporated in the modern theory of
motives? How deep are the apparent analogies between wild ramification in
characteristic $p$ and irregularity? Can one define periods for irregular
connections? If so, do the resulting period matrices have anything to do
with $\epsilon$-factors for $\ell$-adic representations? On the geometric
side, one can ask for a theory of characteristic classes $c_i(E,\nabla)$
for irregular connections such that $c_0$ is the rank of the connection
and $c_1$ is the isomorphism class of the determinant. With these, one
can try to attack the Riemann-Roch problem. 

In this note, we describe
some conjectures and examples concerning what might be called a families
index theorem for irregular connections. Let 
\eq{1.1}{ f:X \to S
}
be a smooth, projective family of curves over a smooth base $S$. Let $\sD$
be an effective relative divisor on $X$. Let $E$ be a vector bundle on
$X$, $\nabla:E \to E\otimes\Omega^1_X(\sD)$ an integrable, absolute
connection with poles along $\sD$. The relative de Rham cohomology
$\R f_*(\Omega^*_{X-\sD}\otimes E)$ inherits a connection (Gau\ss-Manin
connection), and the families index or Riemann-Roch problem is to
describe the isomorphism class of the line bundle with connection
\eq{1.2}{\det \R f_*(\Omega^*_{X-\sD}\otimes E).
} 

Notice this is algebraic de Rham cohomology. Analytically (i.e.
permitting coordinate and basis transformations with essential
singularities on $\sD$), the bundle
$E^{{\rm an}}|_{X-\sD}$ can be transformed (locally on $S$) to have regular
singular points along $\sD$, but the algebraic problem we pose is more
subtle. If $\sD=\emptyset$, or, more generally, if $\nabla$ has regular
singular points, the answer is known:
\eq{1.3}{\det \R f_*(\Omega^*_{X-\sD}\otimes E) = f_*((\det E^\vee,-\det
\nabla)\cdot c_1(\omega_{X/S})) 
} 
(In the case of regular singular points, the $c_1$ has to be taken as a
relative class \cite{ST}, \cite{reg}.) 

In a recent article \cite{irreg}, 
we proved an analogous formula in the case when $E$
was irregular and rank
$1$. For a suitable $\sD$, the relative connection induces an isomorphism
$$\nabla_{X/S,\sD} : E|_\sD \cong E|_\sD\otimes(\omega(\sD)/\omega),
$$
i.e. a trivialization of $\omega(\sD)|_\sD$. The connection pulls back from a rank
$1$ connection $(\sE,\nabla_\sE)$ 
on the relative Picard scheme $\text{Pic}(X,\sD)$,
and the Gau\ss-Manin determinant connection is obtained by evaluating
$(\sE,\nabla_\sE)$ at the priviledged point $(\omega(\sD),\nabla_{X/S,\sD})\in
\text{Pic}(X,\sD)$. 

We want now to consider two sorts of generalizations. First, we
formulate an analogous conjecture for higher rank connections
which are admissible 
in a suitable sense.  
We prove two special cases of this conjecture, computing
the determinant of the Fourier transform of an
arbitrary connection and, up to
$2$-torsion, the determinant of cohomology of
the basic rank $2$ Kloosterman sheaf.

Second, we initiate in the very simplest of cases $E=\sO,\
\nabla(1)=df$, the study 
of periods for irregular connections. Let $m=\deg f+1$. We are led to a
stationary phase integral calculation over the subvariety of 
$\text{Pic}(\P^1,m\cdot\infty)$
corresponding to trivializations of $\omega(m\cdot\infty)$ 
at $\infty$. 
The subtle part of the integral is concentrated at the same point
$$(\omega(m\cdot\infty),\nabla_{X/S,m\cdot\infty})\in
\text{Pic}(\P^1,m\cdot\infty)
$$
mentioned above. 

Our hope is that the geometric methods discussed here will carry over to
the sort of arithmetic questions mentioned above. We expect that the distinguished
point given by the trivialization of $\omega(\sD)$ plays some role in the
calculation of $\epsilon$-factors for rank $1$ $\ell$-adic sheaves and that the
higher rank conjectures have some $\ell$-adic interpretation as well. 

We would like to thank 
Pierre Deligne for sharing his unpublished letters with us.
The basic idea of using the relative Jacobian to study $\epsilon$-factors for rank
$1$ sheaves we learned from him. We also have gotten considerable inspiration from
the works of T. Saito and T. Terasoma cited in the bibliography. The
monograph \cite{Katz} is an excellent reference for Kloosterman sheaves,
and the subject of periods for exponential integrals is discussed briefly
at the end of \cite{kont}.

\section{The Conjecture}

Let $S$ be a smooth scheme over a field $k$ of characteristic $0$. We
consider a smooth family of curves
$f:X
\to S$ and a vector bundle
$E$ with an absolute, integrable connection $\nabla:E \to
E\otimes\Omega^1_X(*D)$. Here $D\subset X$ is a divisor which is smooth
over $S$. We are interested in the determinant of the Gau\ss-Manin
connection
\ml{2.1}{\det \R f_*(\Omega^*_{(X-D)/S}\otimes E) \in
\text{Pic}^\nabla(S):=\H^1(S,\sO_S^\times
\stackrel{d\log}{\longrightarrow} \Omega^1_S) \\
\cong
\Gamma(S,\Omega^1_S/d\log\sO^\times_S). } 
\begin{prop}Let $K=k(S)$ be the function field of $S$. Then the
restriction map $\Pic ^\nabla(S) \to \Pic ^\nabla(\Spec(K))$ is
an injection. 
\end{prop}
\begin{proof} In view of the interpretation of $\text{Pic}^\nabla$ as
$\Gamma(\Omega^1/d\log \sO^\times)$, the proposition follows from the
fact that for a meromorphic function $g$ on $S$, we have $g$ regular at a
point $s\in S$ if and only if $\frac{dg}{g}$ is regular at $s$. 
\end{proof}

Thus, we do not lose information by taking the base $S$ to be the
spectrum of a function field, $S=\Spec(K)$. We shall restrict ourselves
to that case. 

Let 
\eq{2.2}{\sD = \sum_{x\in D(\overline{K})}m_xx
}
be an effective divisor supported on $D$. Suppose that the relative connection has
poles of order bounded by $\sD$, i.e. 
\eq{2.3}{\nabla_{X/S} : E \to E\otimes \omega(\sD),
}
and that $\sD$ is minimal with this property. 
(We frequently write $\omega$ in place of $\Omega_{X/S}$.) We view $\sD$ as an
artinian scheme with sheaf of functions $\sO_\sD$ and 
relative dualizing sheaf $\omega_\sD :=
\omega(\sD)/\omega$. (To simplify, we do not use the notation 
$\omega_{\sD/K}$).  
\begin{prop}\label{prop2.2}Define $E_\sD := E\otimes \sO_\sD$. then $\nabla_{X/S}$
induces a function-linear map
$$\nabla_{X/S,\sD} : E_\sD \to E_\sD\otimes \omega_\sD
$$
\end{prop}
\begin{proof}Straightforward.
\end{proof}

Let $j:X-D \inj X$. There are now two relative
de Rham complexes we might wish to study
\begin{gather*} j_*j^*E \to j_*j^*(E\otimes \omega) \\
E \to E\otimes \omega(\sD).
\end{gather*}
The first is clearly the correct one. 
For example, its relative cohomology carries the Gau\ss-Manin
connection. On the other hand, the second complex is sometimes 
easier to study, as the relative cohomology of the sheaves
involved has finite dimension over $K$.
Since we are primarily interested in the irregular case, that
is when points of $\sD$ have multiplicity $\ge 2$, the following proposition
clarifies the situation. 
\begin{prop}\label{prop2.3} Let notation be as above, and assume every point of $\sD$
has multiplicity $\ge 2$. Then the natural inclusion of complexes
$$\iota:\{E\to E\otimes\omega(\sD)\} \inj \{j_*j^*E \to j_*j^*(E\otimes \omega)\}
$$
is a quasiisomorphism if and only if $\nabla_{X/S,\sD} : E_\sD \to E_\sD\otimes
\omega_\sD$ in proposition \ref{prop2.2} is an isomorphism. 
\end{prop}
\begin{proof} Let $\sD:h=0$ be a local defining equation for $\sD$. Write $\sD =
\sD'+D$ where $D$ is the reduced divisor with support $=\text{supp}(\sD)$.  We claim
first that the map
$\iota$ is a quasiisomorphism if and only if for all $n\ge 1$ the map
$G$ defined by the commutative diagram
$$\begin{CD} E/E(-\sD) @>\nabla_{X/S,\sD} -n\cdot id\otimes \frac{dh}{h} >>
(E(\sD')/E(-D))\otimes
\omega(D) \\ @V\cong V``\cdot h^{-n}" V @V\cong V ``\cdot h^{-n}" V \\
\makebox[1cm][c]{$E(n\sD)/E((n-1)\sD)$} @>G>>
\makebox[1.5cm][c]{$\ \ \ \ \ \ \Big(E(n\sD+\sD')/E((n-1)\sD+\sD')\Big)\! \otimes\!
\omega(D)$}
\end{CD}
$$
given by $h^{-n}e\mapsto h^{-n}\nabla_{X/S,\sD}(e)-nh^{-n-1}e\otimes dh$ is a
quasi-isomorphism. This follows by considering the cokernel of $\iota$
$$j_*j^*E/E \to j_*j^*(E\otimes\omega)/E\otimes\omega(\sD)
$$
and filtering by order of pole. The assertion of the proposition follows because
$nh^{-n-1}e\otimes dh$ has a pole of order strictly smaller than the multiplicity of
$(n+1)\sD$ at every point of
$\sD$. 
\end{proof}

We will consider only the case
\eq{2.4}{\nabla_{X/S,\sD} : E_\sD \cong E_\sD \otimes \omega_\sD.
}

We now consider the sheaf $j_*j^*\Omega^1_X$ of absolute (i.e. relative to $k$)
$1$-forms on $X$ with poles on $D$. Let $\sD=\sum m_xx$ be an effective divisor
supported on $D$ as above, and write $\sD'= \sD-D$. 
\begin{defn} The sheaf $\Omega^p_X\{\sD\}\subset j_*j^*\Omega^p_X$ is defined
locally around a point $x\in D$ with local coordinate $z$ by 
$$\Omega^p_X\{\sD\}_x = \Omega^p_X(\sD')_x +
\Omega^{p-1}_{X,x}\wedge\frac{dz}{z^{m_x}}.
$$
\end{defn}

The graded sheaf $\bigoplus_p\Omega^p_X\{\sD\}$ is stable under the exterior
derivative and independent of the choice of local coordinates at the points of
$D$. One has exact sequences
\eq{2.5}{0 \to f^*\Omega^p_S(\sD') \to \Omega^p_X\{\sD\} \to
\Omega^{p-1}_S(\sD')\otimes \omega(D) \to 0. 
}
\begin{defn}\label{defadm} An integrable absolute connection on $E$ will be called
admissible if there exists a divisor $\sD$ such that $\nabla:E \to
E\otimes\Omega^1_X\{\sD\}$ and such that $\nabla_{X/S,\sD}:E_\sD \cong
E_\sD\otimes\omega_\sD$. 
\end{defn}

\begin{remark} When $E$ has rank $1$, there always exists a $\sD$ such that $\nabla$
is admissible for $\sD$ (\cite{irreg}, lemma 3.1). In higher rank, this need not be
true, even if $\nabla_{X/S,\sD}$ is an isomorphism for some $\sD$. For example, let
$\eta \in \Omega^1_K$ be a closed $1$-form. Let $n\ge 1$ be an integer and let
$c\in k,\ c\neq 0$. The connection matrix 
$$A=\begin{pmatrix}\frac{cdz}{z^m} & \frac{\eta}{z^n} \\ 0 &
\frac{cdz}{z^m}-\frac{ndz}{z}\end{pmatrix}
$$
satisfies $dA+A\wedge A=0$ for all $m, n \in \N$, 
but the resulting integrable connection is not
admissible for $n>m$, although $\nabla_{X/S,(0)}$ is an
isomorphism for $c\neq n$ if $m=1$. Note in this case it is
possible to change basis to get an admissible connection. We don't know what to
expect in general. There do exist connections for which
$\nabla_{X/S,\sD}$ is not an isomorphism for any $\sD$, for example, if
one takes a sum of rank 1 connections with different $m_x$ as above (see
notations \eqref{2.2}) and a local basis adapted to this 
direct sum decomposition. 
\end{remark}

Henceforth, $S=\Spec(K)$ is the spectrum of a function field, and we consider only
integrable, admissible connections
$\nabla:E\to E\otimes \Omega^1_X\{\sD\}$ with $\sD\neq \emptyset$. In
sections $3$ and $4$ we will see many importnat examples (Fourier
transforms, Kloosterman sheaves) of admissible connections. By abuse of
notation, we write
\eq{2.6}{H^*_{DR/S}(E) := \H^*(X,E\to E\otimes\omega(\sD)).
}
We assume (proposition \ref{prop2.3}) 
this group coincides with $H^*_{DR}(E|_{X-D})$.
The isomorphism class of the Gau\ss-Manin connection on the $K$-line $\det
H^*_{DR/S}(E)$ is determined by an element
$$\det H^*_{DR/S}(E)\in\Omega^1_{K/k}/d\log(K^\times)
$$
which we would like to calculate. 

Suppose first that $E$ is a line bundle. Twisting $E$ by
$\sO(\delta)$ for some
divisor $\delta$ supported on the irregular part of the
divisor $D$, we may assume $\deg E=0$. In this case, the result
(the main theorem in
\cite{irreg}) is the following. Since $E$ has rank
$1$,
$\nabla_{X/S,\sD} :E_\sD \cong E_\sD\otimes\omega_\sD$ can be
interpreted as a section
of $\omega_\sD$ which generates this sheaf as an $\sO_\sD$-module.
The exact sequence 
\eq{2.6a}{0 \to \omega \to \omega(\sD) \to \omega_\sD \to 0
}
yields an element $\partial\nabla_{X/S,\sD}\in
H^1(X,\omega)\cong K$ which is known to equal $\deg E=0$. Thus,
we can find some $s\in H^0(X,\omega(\sD))$ 
lifting $\nabla_{X/S,\sD}$. We write $(s)$ for the divisor of
$s$ as a section of 
$\omega(\sD)$ (so $(s)$ is disjoint from $D$). Then the result is
\eq{2.7}{\det H^*_{DR/S}(E) \cong -f_*((s)\cdot E).
}
(When $(s)$ is a disjoint union of $K$-points, the notation on the right simply
means to restrict $E$ with its absolute connection to each of
the points and then 
tensor the resulting $K$-lines with connection together.) Notice that
unlike the classical Riemann-Roch situation (e.g. \eqref{1.3}) the
divisor $(s)$ depends on $(E, \nabla_{X/S})$. 

Another way of thinking about \eqref{2.7} will be important when we consider
periods. It turns out that the connection $(E,\nabla)$ pulls
back from a rank $1$ 
connection $(\sE,\nabla_\sE)$ on the relative Picard scheme
$\text{Pic}(X,\sD)$ whose 
points are isomorphism classes of line bundles on $X$ with
trivializations along 
$\sD$. The pair $(\omega(\sD),\nabla_{X/D,\sD})$ determine a point $t\in
\text{Pic}(X,\sD)(K)$, and \eqref{2.7} is equivalent to
\eq{2.7a}{\det H^*_{DR/S}(E) \cong -(\sE,\nabla_\sE)|_t.
}

Let $\omega_\sD^\times \subset \omega_\sD$ be the subset of elements generating
$\omega_\sD$ as an $\sO_\sD$-module. Let $\partial:\omega_\sD
\to H^1(X,\omega)=K$. 
Define $\tilde B = \omega_\sD^\times \cap \partial^{-1}(0)$. One has a natural
action of $K^\times$, and the quotient $\omega_\sD^\times/K^\times$ is
identified with isomorphism classes of trivializations of
$\omega(\sD)|_\sD$, and hence with a subvariety of
$\text{Pic}(X,\sD)$. One has
\eq{2.7b}{t\in B:=\tilde B/K^\times \subset \omega_\sD^\times/K^\times \subset
\text{Pic}(X,\sD). 
}
The relation between $t, B, \sE$ is the following. $\sE|_B \cong \sO_B$,
so the connection $\nabla_\sD|_B$ is determined by a global
$1$-form $\Xi$. Then 
\eq{2.7c}{\Xi(t) =  0 \in \Omega^1_{B/K}\otimes K(t). 
}
Indeed, $t$ is the unique point on $B$ where the relative 
$1$-form $\Xi/K$ vanishes
(cf \cite{irreg}, Lemma 3.10).

Now suppose the rank of $E$ is $>1$. We will
see when we consider examples in the next section that $\det H^*_{DR/S}(E)$ depends
on more than just the connection on
$\det E$ (remark \ref{rem3.3}). Thus, it is hard to imagine a simple formula like
\eqref{2.7}. Indeed, there is no obvious way other than by taking the
determinant to get rank $1$ connections on $X$ from
$E$. The truly surprising thing is that if we rewrite
\eqref{2.7} algebraically we find a formula which does admit a plausible
generalization. We summarize the results, omitting proofs (which are given in detail
in \cite{irreg}). For each $x_i \in D$, choose a local section $s_i$ of
$\Omega^1_X\{\sD\}$ whose image in $\omega(\sD)$ generates at $x_i$. Write the local
connection matrix in the form
\eq{2.8}{A_i = g_is_i+\frac{\eta_i}{z_i^{m_i-1}}
}
with $z_i$ a local coordinate and $\eta_i\in f^*\Omega^1_S$. Let $s$ be a
meromorphic section of $\omega(\sD)$ which is congruent to the $\{s_i\}$ modulo
$\sD$. Define
\eq{2.9}{c_1(\omega(\sD),\{s_i\}) := (s) \in \Pic(X,\sD).
}
As in \eqref{2.7}, we can define
\eq{2.10}{f_*(c_1(\omega(\sD),\{s_i\})\cdot \det(E,\nabla))\in \Omega^1_K/d\log
K^\times. }
We further define
\eq{2.11}{\{c_1(\omega(\sD)),\nabla\} := f_*\Big(c_1(\omega(\sD),\{s_i\})\cdot
\det(E,\nabla)\Big) -\sum_i \text{res Tr}(dg_ig_i^{-1}A_i).
}
Here res refers to the map 
\eq{2.12}{\Omega^2_X\{\sD\} \to
\Omega^1_K\otimes\omega(\sD)\to\Omega^1_K\otimes\omega_\sD
\overset{\text{transfer}}{\longrightarrow}  \Omega^1_K.
}
\begin{conj}\label{conj2.7}Let $\nabla:E\to E\otimes\Omega^1_X\{\sD\}$ be an
admissible connection as in definition \ref{defadm}. Then
$$\det H^*_{DR/S}(X-D,E) = -\{c_1(\omega(\sD)),\nabla\}\in
\Omega^1_K/d\log(K^\times)\otimes _\Z \Q.
$$
\end{conj}
Our main objective here is to provide evidence for this
conjecture. Of course, one
surprising fact is that the right hand side is independent of
choice of gauge, etc. 
Again, the proof is given in detail in \cite{irreg} but we
reproduce two basic 
lemmas. There are function linear maps
\eq{2.13}{\nabla_{X,\sD}:E_\sD \to
E_\sD\otimes\Big(\Omega^1_X\{\sD\}/\Omega^1_X\Big) ;\quad
\nabla_{X/S,\sD}:E_\sD \to 
E_\sD\otimes\omega_\sD, 
}
and it makes sense to consider the commutator
$$[\nabla_{X,\sD},\nabla_{X/S,\sD}]:E_\sD \to E_\sD\otimes
\Big(\Omega^1_X\{\sD\}/\Omega^1_X\Big)\otimes\omega_\sD.
$$
\begin{lem}\label{comlem} $[\nabla_{X,\sD},\nabla_{X/S,\sD}]=0.$ \newline
\end{lem}
\begin{proof}With notation as in \eqref{2.8}, we take $s_i =
\frac{dz_i}{z_i^{m_i}}$, where $z_i$ is a local coordinate. Integrability implies
\eq{2.14}{dA_i =
dg_i\wedge\frac{dz_i}{z_i^{m_i}}+d\Big(\frac{\eta_i}{z_i^{m_i-1}}\Big)=A_i^2 =
[\eta_i,g_i]\frac{dz_i}{z_i^{2m_i-1}} + \epsilon 
}
with $\epsilon \in \Omega^2_K\otimes K(X)$. 
Multiplying through by $z_i^{m_i}$, we conclude
that $[\frac{\eta_i}{z_i^{m_i-1}},g_i]$ is regular on $D$, which is equivalent to
the assertion of the lemma. 
\end{proof}
The other lemma which will be useful in evaluating the right hand term in
\eqref{2.11} is
\begin{lem}\label{lem2.9}We consider the situation from
\eqref{2.8} and \eqref{2.11} 
at a fixed
$x_i \in D$. For simplicity, we drop the $i$ from the notation.
Assume $ds=0$. Then 
$$\rt (dgg^{-1}A) = \rt (dgg^{-1}\frac{\eta}{z^{m-1}}). 
$$
\end{lem}
\begin{proof}We must show $\text{res Tr}(dgs)=0$. Using \eqref{2.14} and
$\Tr[g,\eta]=0$ we reduce to showing $0=\rt (d(\eta z^{1-m})) \in
\Omega^1_K$. Since $\eta\in f^*\Omega^1_K$, we may do the computation formally
locally and replace $d$ by $d_z$. The desired vanishing follows because an
exact form has no residues. 
\end{proof}

 With lemma \ref{comlem}, we can formulate the conjecture in a more invariant
way in terms of an AD-cocycle on $X$. Recall \cite{E}
\eq{2.15}{AD^2(X) := \H^2(X,\sK_2
\stackrel{d\log}{\longrightarrow} \Omega^2_X)\cong
H^1(X,\Omega^2_X/d\log\sK_2) 
}
The AD-groups are the cones of cycle maps from Chow groups to Hodge cohomology,
and as such they carry classes for bundles with connections. There is a general
trace formalism for the AD-groups, but in this simple case the reader can easily
deduce from the right hand isomorphism in \eqref{2.15} a trace map
\eq{2.16}{ f_* : AD^2(X) \to AD^1(S) = \Omega^1_K/d\log K^\times.
}
When the connection $\nabla$ has no poles (or more generally, when it has regular
singular points) it is possible to define a class
\eq{2.17}{\epsilon = c_1(\omega)\cdot c_1(E,\nabla)\in AD^2(X)
}
with $f_*(\epsilon) = [\det H^1_{DR/S}(X,E)]$. Remarkably, though it no longer has
the product description \eqref{2.17}, one can associate such a class to any
admissible connection. 

Fix the divisor $\sD$ and consider tuples $\{E,\nabla,\sL,\mu\}$
where $(E,\nabla)$ is
an admissible, absolute connection, $\sL$ is a line bundle on
$X$, and $\mu:E_\sD
\cong E_\sD\otimes \sL_\sD$. We require
\eq{2.18}{0=[\mu,\nabla_{X,\sD}]:E_\sD \to
E_\sD\otimes\sL_\sD\otimes\Big(\Omega^1_X\{\sD\}/\Omega^1_X\Big). }
Of course, the example we have in mind, using lemma \ref{comlem}, is
\eq{2.19}{\{E,\nabla\} := \{E,\nabla,\omega(\sD),\nabla_{X/S,\sD}\}.
}
To such a tuple satisfying \eqref{2.18}, we associate a class
$\epsilon(E,\nabla,\sL,\mu)\in AD^2(X)$ as follows. Choose cochains $c_{ij}\in
GL(r,\sO_X)$ for $E$, $\lambda_{ij}\in \sO_X^\times$ for $\sL$, $\mu_i \in
GL(r,\sO_\sD)$ for $\mu$, and 
$\omega_i \in M(r\times r,\Omega^1_X\{\sD\})$ for
$\nabla$. Choose local liftings $\tilde{\mu_i}\in GL(r,\sO_X)$ for the $\mu_i$.

\begin{prop} \label{coc} The Cech hypercochain
$$\Big(\{\lambda_{ij},\det(c_{jk})\},d\log \lambda_{ij}\wedge
\Tr(\omega_j),\Tr(-d\tilde{\mu_i}\tilde{\mu_i}^{-1}\wedge \omega_i)\Big)
$$
represents a class 
$$\epsilon(E,\nabla,\sL,\mu)\in \H^2(X,\sK_2 \to \Omega^2_X\{\sD\} \to
\Omega^2_X\{\sD\}/\Omega^2_X)
\cong AD^2(X).
$$
This class is well defined independent of the various choices. Writing
$\epsilon(E,\nabla) = \epsilon(E,\nabla,\omega(\sD),\nabla_{X/S,\sD})$, we have
$$f_*\epsilon(E,\nabla) = \{c_1(\omega(\sD)),\nabla\} 
$$
where the right hand side is defined in \eqref{2.11}.
\end{prop}
\begin{proof} Again the proof is given in detail in \cite{irreg}
and we omit it. 
\end{proof}

As a consequence, we can restate the main conjecture:
\begin{conj}\label{conj2.11} Let $\nabla$ be an integrable,
admissible, absolute 
connection as above. Then
$$\det H^*_{DR/S}(E,\nabla) = -f_*\epsilon(E,\nabla).
$$
\end{conj}
 
To finish this section, we would like to show that behind the
quite technical cocyle written in proposition \ref{coc}, there
is an algebraic group playing a role similar to $\Pic(X, \sD)$
in the rank 1 case. Let $G$ be the algebraic group whose $K$-points are 
isomorphism classes $(\sL,\mu)$, where $\sL$ is an invertible
sheaf, and $\mu: E_{\sD} \to E_{\sD}\otimes \sL_{\sD}$ is an
isomorphism commuting with $\nabla_{X,\sD}$. It is endowed with
a surjective map $q: G\to \Pic(X)$.
As noted, $G$
contains the special point $(\omega(\sD),
\nabla_{X/S, \sD})$. 
The cocyle of proposition \ref{coc}
defines a class in $\H^2(X\times_K G, \sK_2 \to
\Omega^2_{X\times G}\{\sD\times G\} \to     
\Omega^2_{X\times G}\{\sD\times G\}/\Omega^2_{X\times G})$.
Taking its trace \eqref{2.16}, one obtains a class in
$(\sL(E),\nabla(E)) \in AD^1(G)$,
that is a rank one connection
on $G$. Then $f_*\epsilon(E,\nabla)$ is simply the restriction
of $(\sL(E),\nabla(E))$ to the special point 
$(\omega(\sD), \nabla_{X/S, \sD})$.

Now we want to show that this special point, as in the rank 1
case, has a very special meaning.   By analogy with \eqref{2.7b}, we
define
\begin{gather}
\tilde B= \Big(\Ker ({\rm Hom}(E_\sD, E_\sD \otimes
\omega_\sD)\xrightarrow{\rt} K)\Big)\cap \text{Isom}(E_\sD, E_\sD \otimes
\omega_\sD) \\ 
B = \tilde B/K^\times \subset G.
\end{gather}
We observe that lemma \ref{lem2.9} shows that $g\in B$. 
Choosing  
a local trivialization $\omega_\sD\cong
\sO_\sD \frac{dz}{z^m}$ and a local trivialization of $E_\sD$,
we write
$\nabla_{X/S}|_{\sD}$ as a matrix $g\frac{dz}{z^m}$, with $g\in
GL_r(\sO_{\sD})$. $\Theta$ is then identified with a translation invariant
form on the restriction of scalars $\text{Res}_{\sD/K}GL_r$
\begin{gather}
\Theta:= 
(\sL(E),\nabla_{G/K}(E))= \rt (d\mu \mu^{-1} g \frac{dz}{z^m}).
\end{gather}
The assumption that $g\in B$ implies that $\Theta$ descends to an
invariant form on $\text{Res}_{\sD/K}GL_r/\G_m\supset G_0$, where $G_0 :=
\{(\sO,\mu)\in G\}$. By invariance, it gives rise to a form on the
$G_0$-torsor $G_{\omega(\sD)} := \{(\omega(\sD),\mu)\in G\}$. We have
$$ B \subset G_{\omega(\sD)} \subset G.
$$
Let $S\subset G_0$ be the subgroup of points stabilizing $B$. 
\begin{prop}$\Theta|_B$ vanishes at a point $t\in B$ if and only if $t$
lies in the orbit $g\cdot S$. 
\end{prop}
\begin{proof} Write $th=g$. Write the universal element in
$\text{Res}_{\sD/K}GL_r$ as a matrix
$X=\sum_{k=0}^{m-1}(X_{ij}^{(k)})_{ij}z^k$. The assertion that $\Theta|_B$
vanishes in the fibre at $t$ means 
$$\text{res Tr}(\sum_k(d(X_{ij}^{(k)})_{ij}z^k)h\frac{dz}{z^m})(t) =
a(\sum_idX_{ii}^{(m-1)})(t) 
$$
for some $a\in K$. Note this is an identity of the form 
\eq{elts}{0 = \sum c_{ij}^{(k)}dX_{ij}^{(k)}\in \Omega^1_{G/K}\otimes
K(t);\quad c_{ij}^{(k)}\in K. 
}
We first claim that in fact this identity holds already in
$\Omega^1_{G/K}$. To see this, write $\sG=\text{Res}_{\sD/K}GL_r$. Note
$\Omega^1_\sG$ is a free module on generators $dX_{ij}^{(k)}$.  Also,
$G\subset \sG$ is defined by the equations
\begin{gather}
[\sum_{k=0}^{m-1} (X_{ij}^{(k)})_{ij}z^k,\sum_{k=0}^{m-1}
( g_{ij}^{(k)})_{ij}z^k]=0 \notag \\
[\sum_{k=0}^{m-1} (X_{ij}^{(k)})_{ij}z^k,\sum_{k=0}^{m-2}
( \eta_{ij}^{(k)})_{ij}z^k]=0,\notag
\end{gather}
which are of the form
$$\sum_{i,j,k} b_{ijp}^{(k)}X_{ij}^{(k)}=0,\quad p=1,2,\ldots,M;\quad 
b_{ijp}^{(k)}\in K,
$$
that is are linear equations in the $X_{ij}^{(k)}$ with
$K$-coefficients. 
Thus, we have an exact sequence
\eq{exseq}{0 \to N^\vee \to \Omega^1_{\sG/K}\otimes\sO_G \to \Omega^1_{G/K} \to
0, 
}
where $N^\vee$ is generated by $K$-linear combinations of the $dX_{ij}^{(k)}$.
We have, therefore, a reduction of structure of the sequence \eqref{exseq}
from $\sO_G$ to $K$, and therefore $\Omega^1_{G/K} \cong
\Omega^1_0\otimes_K\sO_G$, where $\Omega^1_0 \subset \Omega^1_{G/K}$ is
the $K$-span of the $dX_{ij}^{(k)}$. Hence, any $K$-linear identity among
the $dX_{ij}^{(k)}$ which holds at a point on $G$ holds everywhere on
$G$. 
As a consequence, we can integrate to an identity
\eq{main1}{\text{res Tr}(Xh\frac{dz}{z^m}) =
a(\sum_iX_{ii}^{(m-1)}) + \kappa
}
with $\kappa\in K$. If we specialize $X\to t$ we find
\eq{main2}{0 = \text{res Tr}(g\frac{dz}{z^m}) = a\cdot \text{res
Tr}(t\frac{dz}{z^m}) +\kappa = \kappa
}
We conclude from \eqref{main1} and \eqref{main2} that $h\in S$.
\end{proof}

\section{The Fourier Transform}

In this section we calculate the Gau\ss-Manin determinant line for the Fourier
transform of a connection on $\P^1-D$ and show that it satisfies the conjecture
\ref{conj2.7}. Let $\sD=\sum m_\alpha\alpha$ be an effective $k$-divisor on
$\P^1_k$. Let $\sE=\oplus_r \sO$ be a rank $r$ free bundle on $\P^1_k$, and let
$\Psi:\sE \to \sE\otimes \omega(\sD)$ be a $k$-connection on
$\sE$. Let $(\sL,\Xi)$
denote the rank $1$ connection on $\P^1\times\P^1$ with poles on
$\{0,\infty\}\times 
\{0,\infty\}$ given by 
$\sL = \sO_{\P^1\times \P^1}$ and
$\Xi(1) = d\big(\frac{z}{t}\big)$. Here $z,t$ are the
coordinates on the two copies of the projective line. Let
$K=k(t)$. We have a diagram 
\eq{3.1}{\begin{array}{ccc} (\P^1_z-\sD)\times \P^1_t & \hookleftarrow &
(\P^1_z-\sD)\times
\Spec(K)
\\
\rule{0cm}{.5cm}\makebox[.4cm][l]{$\downarrow p_1$} &&
\makebox[.4cm][l]{$\downarrow 
p_2$} \\
\rule{0cm}{.5cm}\P^1_z-\sD && \Spec(K)
\end{array}
}
The Gau\ss-Manin determinant of the Fourier transform is given
at the generic point 
by
\ml{3.2}{\det H^*_{DR/K}\Big((\P^1_z-\sD)_K,\
p_1^*(\sE,\Psi)\otimes (\sL,\Xi)\Big) 
\\
= \det H^*_{DR/K}\Big((\P^1_z-\sD)_K,(E,\nabla)\Big)
}
with $E:=p_1^*\sE\otimes \sL|_{(\P^1_z-\sD)_K}$ and
$\nabla=\Psi\otimes 1 + 1\otimes 
\Xi$. We have the following easy
\begin{remark} Write
\eq{3.3}{\Psi  = \sum_\alpha \sum_{i=1}^{m_\alpha}\frac{g^\alpha_i
dz}{(z-\alpha)^i} +d(g^\infty_1z+\ldots+g^\infty_{m_\infty-1}z^{m_\infty-1}) 
}
where $g^\alpha_i\in M(r\times r,k)$.
Then 
\eq{3.4}{\nabla=\Psi+\frac{dz}{t}-\frac{zdt}{t^2}
}
is admissible if and only if either
\begin{enumerate}\item[(i)] $m_\infty \le 2$ and
$g^\alpha_{m_\alpha}$ is invertible 
for all $\alpha\neq \infty$,  or
\item[(ii)] $m_\infty\ge 3$, $g^\alpha_{m_\alpha}$ is invertible
for all $\alpha\neq \infty$, and $g^\infty_{m_\infty-1}$ is invertible. 
\end{enumerate}
\end{remark}
\begin{thm}The connection $(E,\nabla)$ satisfies conjecture \ref{conj2.7}.
\end{thm}
\begin{proof} We first consider the case when $\Psi$ has a pole
of order $\le 1$ at 
infinity, so the $g^\infty_i=0$ in \eqref{3.3}. A basis for 
$$
H^0(\P^1_K, E\otimes\omega(\sum m_\alpha\alpha+2\infty))
$$
is given by
\eq{3.5}{e_j\otimes dz;\quad e_j\otimes \frac{dz}{(z-\al)^i},\
1\le i\le m_\al,\ 
1\le j\le r. 
}
$H^0_{DR/K} = (0)$ and $H^1_{DR/K} = \text{coker}(H^0(E) \to
H^0(E\otimes\omega(\sum m_\al\al+2\infty)))$ has basis
\eq{3.6}{e_j\otimes \frac{dz}{(z-\al)^i};\quad 1\le i\le m_\al,\ 1\le j\le r.
}
To compute the Gau\ss-Manin connection, we consider the diagram (here $\sD=\sum
m_\al\al +2\infty$ and $\sD' = \sD-D =\sum(m_\al-1)\al+\infty$)
\begin{tiny}
\eq{3.7}{\minCDarrowwidth.4cm\begin{CD}
@. @. H^0(E) @= H^0(E) \\
@. @. @VV\nabla_XV @VV\nabla_{X/S}V \\
0 @>>> H^0(E(\sD'))\otimes\Omega^1_K  @>>> H^0(E\otimes\Omega^1_{\P^1}\{\sD\}) @>a>>
H^0(E\otimes\omega(\sD)) @>>> 0 \\
@. @VV\nabla_{X/S}\otimes 1 V @VV\nabla_X V\\
@. H^0(E(\sD')\otimes\omega(\sD))\otimes\Omega^1_K @>\cong >> 
H^0(E\otimes\Omega^2_{\P^1}\{\sD+\sD'\})
\end{CD}
}
\end{tiny}
One deduces from this diagram the Gau\ss-Manin connection 
\eq{3.8}{H^1_{DR/K}(E)\cong \text{coker}(\nabla_{X/S})
\xrightarrow{\nabla_{GM}} 
H^1_{DR/K}(E)\otimes\Omega^1_K;\quad w\mapsto \nabla_X(a^{-1}(w)).  
}
We may choose $a^{-1}(e_j\otimes\frac{dz}{(z-\al)^i})=
e_j\otimes\frac{dz}{(z-\al)^i}$, so by \eqref{3.4}
\eq{3.9}{\nabla_{GM}\Big(e_j\otimes\frac{dz}{(z-\al)^i}\Big) =
\nabla_X\Big(e_j\otimes\frac{dz}{(z-\al)^i}\Big) = e_j\otimes\frac{zdz\wedge
dt}{(z-\al)^it^2}  
}
In $H^1_{DR/K}\cong \text{coker}(\nabla_{X/S})$ we have the identity
\eq{3.10}{e_j\otimes dz = -t\Psi e_j.
}
We conclude
\ml{3.11}{\nabla_{GM}\Big(e_j\otimes\frac{dz}{(z-\al)^i}\Big) \\
= \begin{cases} \Big(e_j\otimes\frac{dz}{(z-\al)^{i-1}}+\al e_j\otimes
\frac{dz}{(z-\al)^{i}}\Big)\wedge \frac{dt}{t^2} & 2\le i\le m_\al \\
-t\Psi e_j +\al e_j\otimes
\frac{dz}{z-\al}\Big)\wedge \frac{dt}{t^2} & i=1.
\end{cases}
}
In particular, the determinant connection, which is given by
$\Tr \nabla_{GM}$, can 
now be calculated:
\eq{3.12}{\Tr\nabla_{GM} = \sum_\al \frac{rm_\al \al dt}{t^2}-\Tr\sum_\al
\frac{g^\al_1 dt}{t}. 
}
We compare this with the conjectured value which is the negative of \eqref{2.11}.
Define
\eq{3.13}{ F(z) := \sum_\al \frac{1}{(z-\al)^{m_\al}} - 1 =
\frac{G(z)}{(z-\al)^{m_\al}} ;\quad s:= F(z)dz.
}
One has
\eq{3.14}{c_1(\omega(\sD),\Big\{\frac{dz}{(z-\al)^{m_\al}},dz\Big\}) = (G),
}
the divisor of zeroes of $G$. We need to compute $(\det
E,\det\nabla)|_{(G)}$. We 
have
\ml{3.15}{G(z) = \sum_\al \prod_{\beta\neq \al}(z-\beta)^{m_\beta}
-\prod_\al(z-\al)^{m_\al} \\
= -z^{\sum m_\al}+\Big(\sum m_\al\al+\#\{\al\,|\,m_\al=1\}\Big)z^{(\sum
m_\al)-1}+\ldots } 
Note
that the coefficients of $G$ do not involve $t$, so the $dz$
part of the connection 
dies on $(G)$ and we get 
\ml{3.16}{\Tr \nabla |_{(G)} = \frac{-rzdt}{t^2}|_{(G)} =
-\frac{rdt}{t^2}\sum_{\substack{\beta \\ G(\beta)=0}} \beta \\
=-\frac{rdt}{t^2}\Big(\sum m_\al\al+\#\{\al\,|\,m_\al=1\}\Big) 
}
It remains to evaluate the correction terms $\text{res
Tr}(dgg^{-1}A)$ occurring in
\eqref{2.11}. In the notation of \eqref{2.8}, $\frac{\eta}{z^{m-1}} =
\frac{-zdt}{t^2}$, and by lemma \ref{lem2.9} we have $\text{res
Tr}(dgg^{-1}A)=-\text{res Tr}(dgg^{-1}\frac{zdt}{t^2})$. Clearly, the only
contribution comes at $z=\infty$. Take $u=z^{-1}$. At $\infty$ the connection is
\eq{3.17}{A=-\Big(\sum_\al\sum_i\frac{g^\al_i
u^i}{(1-u\al)^i}+\frac{1}{t}\Big)\frac{du}{u^2}-\frac{dt}{ut^2}. 
}

We rewrite this in the form $A=gs+\frac{\eta}{u}$ as in \eqref{2.8} with $s$ as in
\eqref{3.13} and $\eta = -\frac{dt}{t^2}$. We find
\eq{3.18}{ g = \frac{\sum_{i,\al}\frac{g^\al_iu^i}{(1-u\al)^i}+t^{-1}}{\sum_\al
\frac{u^{m_\al}}{(1-u\al)^{m_\al}}-1}  = \frac{\kappa}{v},
}
(defining $\kappa$ and $v$ to be the numerator and denominator, respectively.) Then
\ml{3.19}{\text{res Tr}(dgg^{-1}A) = -\text{res Tr}(dgg^{-1}\frac{dt}{ut^2}) \\
=
\Big(-\text{res Tr}(d\kappa\kappa^{-1}u^{-1})+r\cdot\text{res
Tr}(dvv^{-1}u^{-1})\Big)\frac{dt}{t^2} \\
=\Big(-t\sum\Tr(g^\al_1)-r\#\{\al\,|\,m_\al=1\}\Big)\frac{dt}{t^2}
} 
Combining \eqref{3.19}, \eqref{3.16}, and \eqref{3.12} we conclude
\eq{3.20}{\Tr \nabla_{GM} = -\Big(\Tr\nabla|_{(G)} - \text{res Tr}(dgg^{-1}A)\Big),
}
which is the desired formula.

We turn now to the case where $\Psi$ has a pole of order $\ge 2$ at infinity. We
write
\begin{gather}\label{3.21} \Psi =
\sum_\al\sum_{i=1}^{m_\al}\frac{g^\al_idz}{(z-\al)^i} +g^\infty dz;\quad g^\infty =
g^\infty_2+\ldots+g^\infty_{m_\infty}z^{m_\infty-2} \\
\nabla = \Psi + \frac{dz}{t}-\frac{zdt}{t^2}.
\end{gather}
A basis for $\Gamma(\P^1,\omega(\sum m_\al \al +m_\infty\infty))$ is given by 
\eq{3.23}{e_j\otimes \frac{dz}{(z-\al)^{i}};\ \ 1\le i\le m_\al;\quad
e_j\otimes z^idz;\ 0\le i\le m_\infty -2. }
A basis for the Gau\ss-Manin bundle is given by omitting $e_j\otimes
z^{m_\infty-2}dz$. As in \eqref{3.7}-\eqref{3.9}, the Gau\ss-Manin connection is
\eq{3.24}{w\mapsto zw\wedge\frac{dt}{t^2}.
} 
To compute the trace, note that in $H^1_{DR/K}$, we have if $m_\infty \ge 3$
\ml{3.25}{e_j\otimes z^{m_\infty-2}dz \\
= -(g^\infty_{m_\infty})^{-1}
\Big(\sum_{i,\al}g^\al_i(e_j)\otimes\frac{dz}{(z-\al)^i} + \\
(g^\infty_2(e_j)+\ldots+g^\infty_{m_\infty-1}(e_j)z^{m_\infty-3}+t^{-1}e_j)\otimes
dz\Big).  
}
If $m_\infty=2$,
\eq{3.26}{e_j\otimes dz 
= -(g^\infty_{2}+t^{-1})^{-1}
\Big(\sum_{i,\al}g^\al_i(e_j)\otimes\frac{dz}{(z-\al)^i} \Big).  
}
It follows that
\eq{3.27}{\Tr\nabla_{GM} = \begin{cases} \rule{0cm}{.5cm}\Big(\sum_\al rm_\al\al
-\Tr\big((g^\infty_{m_\infty})^{-1}g^\infty_{m_\infty-1}\big)\Big)\frac{dt}{t^2}  &
m_\infty\ge 3
\\
\rule{0cm}{.7cm}\Big(\sum_\al rm_\al\al-\sum_\al
\Tr\big((g^\infty_{2}+t^{-1})^{-1}g^\al_1\big)\Big)\frac{dt}{t^2} & m_\infty=2
\end{cases}
}

To compute the right hand side in conjecture \ref{conj2.7}, we take as trivializing
section
\eq{3.28}{s= \sum_\al \frac{dz}{(z-\al)^{m_\al}} - z^{m_\infty-2}dz =
\frac{G(z)dz}{\prod_\al(z-\al)^{m_\al}}  
}
where 
\eq{3.29}{G(z) = \sum_\al\prod_{\beta\neq \al}
(z-\beta)^{m_\beta} - z^{m_\infty-2}\prod_\al(z-\al)^{m_\al} 
}
We have $(s) = (G)$, the divisor of zeroes of the polynomial $G$. Again, $G$ does
not involve $t$, so if $G=-\prod(z-a_k)$, we have 
\eq{3.30}{\Tr\nabla|_{z=a_p} = -\frac{ra_pdt}{t^2}
}
Thus
\eq{3.31}{\Tr\nabla|_{(s)} = \begin{cases}-\frac{(\sum
rm_\al\al+r\cdot\#\{\al\,|\,m_\al=1\}) dt}{t^2}, & m_\infty=2 \\
-\frac{\sum
rm_\al\al dt}{t^2}, & m_\infty\ge 3.
\end{cases}
}
Finally we have to deal with the correction term $\text{res
Tr}(dgg^{-1}(-\frac{zdt}{t^2}))$. There is no contribution except at $\infty$. We
put $u=z^{-1}$ as before, and $g=\frac{\kappa}{v}$ with 
\begin{gather}\kappa =
\sum_\al\sum_{i=1}^{m_\al}\frac{u^{m_\infty-2+i}g^\al_i}{(1-\al u)^i}
+(g^\infty_2+t^{-1})u^{m_\infty-2}+\ldots+g^\infty_{m_\infty} \\ v=\sum_\al
\frac{u^{m_\infty-2+m_\al}}{(1-\al u)^{m_\al}} - 1
\end{gather}
Assume first $m_\infty = 2$. Then
\ml{3.34}{\text{res Tr}(dgg^{-1}A) = -\text{res Tr}(dgg^{-1}\frac{dt}{ut^2}) \\
=
-\text{res Tr}(d_z\kappa\kappa^{-1}u^{-1})\frac{dt}{t^2} +r\cdot\text{res
Tr}(d_zvv^{-1}u^{-1})\frac{dt}{t^2}  \\
= -\Tr \Big((g^\infty_2+t^{-1})^{-1}\sum g^\al_1\Big)\frac{dt}{t^2}-r\cdot
\#\{\al\,|\,m_\al=1\}\frac{dt}{t^2} 
}
In the case $m_\infty \ge 3$ we find
\eq{3.35}{\text{res Tr}(dgg^{-1}A)
=-\Tr\Big((g^\infty_{m_\infty})^{-1}g^\infty_{m_\infty-1}\Big) \frac{dt}{t^2} 
}
The theorem follows by comparing \eqref{3.27}, \eqref{3.31}, \eqref{3.34}, and
\eqref{3.35}. 
\end{proof}
\begin{remark}\label{rem3.3} The presence of nonlinear terms in the $g^\al_i$ in
\eqref{3.27} means that the connection on $\det H^*_{DR/K}(E)$ is not determined by
$\det E$ alone. 
\end{remark}
\section{Kloosterman Sheaves}

In this section we show that the main conjecture holds at least
up to $2$-torsion for 
the basic rank
$2$ Kloosterman sheaf \cite{Katz}.  
The base field $k$ is $\C$.
Fix $a,b\in K^\times$ (in fact one can work over $K=\C[a,b,a^{-1},b^{-1}]$).
Let $\al, \be\in 
\C-\Z$ and assume also $\al -\be \in \C-\Z$. Consider two connections
$(\sL_i,\nabla_i)$ on the trivial bundle on
$\G_m$ given by $1\mapsto \al d\log(t)+d(at)$ and $1\mapsto \be
d\log(u)+d(bu)$ where 
$t,u$ are the standard parameters on two copies of $\G_m$. Let
$X:= \G_m\times \G_m$ 
and consider the exterior tensor product connection on $X$ 
\eq{4.1}{\sL := \sL_1\boxtimes\sL_2 = (\sO_X,\nabla);\quad
\nabla(1)=\al d\log(t)  + \be d\log(u) + d(at + bu).
}
Note all the above are integrable, absolute connections.
\begin{prop}\begin{enumerate}\item 
$$H^i({DR}(\sL_i/K) \cong \begin{cases}K & i=1 \\
0 & i\not = 1
\end{cases}
$$
\item 
$$H^p_{DR}(\sL/K) \cong \begin{cases} H^1_{DR}(\sL_1/K)\otimes H^1_{DR}(\sL_2/K) & p=2 \\
0 & p\not = 2
\end{cases}
$$
\end{enumerate}
\end{prop}
\begin{proof} The de Rham complex (of global sections on $X$) for $(\sL,\nabla)$ is the
tensor product of the corresponding complexes for the $\nabla_i$, so (2) follows from (1).
For (1) we have e.g. the complex of global sections $\sO
\stackrel{\nabla_{1,K}}{\longrightarrow}\omega,\ 1\mapsto \al d\log(t)+adt$. In $H^1_{DR}$
this gives for all $n\in \Z$
$$ at^ndt \equiv -(\al +n)t^{n-1}dt.
$$
Assertion (1) follows easily.
\end{proof}

We now compute Gau\ss-Manin. We will (abusively) use sheaf
notation when we mean 
global sections over
$\G_m$ or $X$. Also we write $\nabla$ for either one of the
$\nabla_i$ or the exterior 
tensor connection. $\nabla_K$ is the corresponding relative
connection. One has the 
diagram
$$\begin{CD}@. \sO @= \sO \\
@. @VV\nabla V @VV\nabla_K V @. \\
\sO\otimes \Omega^1_K @>>> \Omega^1
@>\stackrel{\sigma}{\longleftarrow}>> \omega \\ 
@VV\nabla_K\otimes 1 V @VV\nabla V @. \\
\omega\otimes \Omega^1_K @> \iota >\cong > \Omega^2/F^2.
\end{CD}
$$
Here $\sigma$ is the obvious function linear map (e.g. for
$\nabla=\nabla_1,\ \sigma(t^ndt)=t^ndt$), and $F^2\subset
\Omega^2$ is the subgroup of 2 forms coming from the base.
This leads to the Gau{\ss}-Manin diagram
$$\begin{CD}\sO @>\nabla - \sigma\nabla_K >> \sO\otimes\Omega^1_K \\
@VV \nabla_K V @VV \nabla_K \otimes 1 V \\
\omega @> -\iota^{-1}\nabla\sigma >> \omega\otimes\Omega^1_K.
\end{CD}
$$
For example when $\nabla=\nabla_1$ we get on $H^1_{DR}$
$$\nabla_{GM}(t^ndt) = -\iota^{-1}(\nabla(1)\wedge t^ndt) =
-\iota^{-1}(tda\wedge t^ndt) = 
t^{n+1}dt\otimes da.
$$
Since $tdt \equiv \frac{-(\al +1)}{a}dt$, we get
$$\nabla_{GM}(dt) = -(\al +1)dt\otimes d\log(a) \equiv -\al
dt\otimes d\log(a) \!\!\mod 
d\log(K^\times).
$$
On the (rank $1$) tensor product connection $H^2_{DR}(\sL) =
H^1_{DR}(\sL_1)\otimes 
H^1_{DR}(\sL_2)$, the Gau\ss-Manin determinant connection is therefore
\begin{equation}\label{1}-\al d\log(a) - \be d\log(b). 
\end{equation}

Note that we computed this determinant here by hand, but we
could have as well applied  directly theorem 4.6 of
\cite{irreg}: the determinant of
$H^1_{DR}(\sL_1)$ is just the restriction of $\nabla_1$ to the
divisor of $\P^1$ defined by the trivializing section $\alpha
d\log(t) + adt$ of $\omega(0 + 2\infty)$,
that is by $\frac{a}{t}
+\alpha=0$. Thus the determinant is $-\al d\log(a)\in
\Omega^1_K/d\log K^\times$,  and similarly for
$H^1_{DR}(\sL_2)$, the determinant is $- \be d\log(b)\in
\Omega^1_K/d\log K^\times$.

The idea now is to recalculate that determinant connection using
the Leray spectral sequence for the map $\pi:X 
\to \G_m,\ \pi(t,u)=tu$. Write $v$ for the coordinate on the
base, so $\pi^*(v)=tu$. 

\begin{prop}We have
$$R^i\pi_{*,DR}(\sL) = \begin{cases} 0 & i\not = 1 \\
\text{rank $2$ bundle on } \G_m & i=1
\end{cases}
$$
\end{prop}
\begin{proof}Let $\nabla_\pi$ be the relative connection on $\sL$ with respect to the map
$\pi$, and take $t$ to be the fibre coordinate for $\pi$. Write $u=\frac{v}{t}$. Then
$$\nabla_\pi(1) = \al d\log(t)+ \be d\log(u) + adt + bdu = (\al - \be)d\log(t)+adt -
bv\frac{dt}{t^2}
$$
so in $R^1\pi_{*,DR}$ we have
$$0\equiv \nabla_\pi(t^n) = (\al - \be+n)t^{n-1}dt+at^ndt -
bvt^{n-2}dt.
$$
It follows that $R^1\pi_{*,DR}(\sL)$ has rank $2$ (generated e.g. by $dt$ and
$\frac{dt}{t}$), and the other $R^i=(0)$ as claimed.
\end{proof}

Define
$$\sE := R^1\pi_{*,DR}(\sL),\ \nabla=\nabla_{GM}:\sE \to
\sE\otimes\Omega^1_{\G_m}. 
$$
\begin{thm}\label{thm4.3}The Gau\ss-Manin connection on $H^*_{DR/K}(\G_m,\sE)$
satisfies conjecture \ref{conj2.7} up to $2$-torsion. 
\end{thm}
\begin{remark}In fact, we will see that $\nabla$ is not
admissible in the sense of 
definition \ref{defadm}, but its inverse image via a degree 2
covering is. Since the new determinant of de Rham cohomology
obtained in this way is twice the old one, 
we lose control of the $2$-torsion. We do not know
whether the conjecture 
holds exactly in this case or not. 
\end{remark}
\begin{proof} We can now calculate the connection $\nabla :=
\nabla_{GM}$ on $\sE$ 
just as before. We have the Gaus{\ss}-Manin 
diagram ($\sigma(t^ndt) = t^ndt$).
$$\begin{CD}\sO_X @>\nabla - \sigma\nabla_\pi >> \sO_X\otimes\Omega^1_{\G_m} \\
@VV \nabla_\pi V @VV \nabla_\pi \otimes 1 V \\
\omega_\pi @> -\iota^{-1}\nabla\sigma >> \omega_\pi\otimes\Omega^1_{\G_m}.
\end{CD}
$$
Here 
$$\nabla(1) = (\al -\be)\frac{dt}{t} + \be \frac{dv}{v} +adt+tda
+t^{-1}d(bv) +bvd(t^{-1}). 
$$
We get
\begin{gather*}\nabla_{GM}(\frac{dt}{t}) =
-\iota(\nabla(1)\wedge\frac{dt}{t}) = 
\frac{dt}{t}\otimes \be \frac{dv}{v} + dt\otimes da +
\frac{dt}{t^2}\otimes d(bv) \\ 
\nabla_{GM}(dt) = -\iota^{-1}((\be \frac{dv}{v} + tda +
t^{-1}d(bv))\wedge dt) = 
\\ dt\otimes
\be\frac{dv}{v} + tdt\otimes da + \frac{dt}{t}\otimes d(bv)
\end{gather*}
We can now substitute
\begin{gather*}\frac{dt}{t^2} \equiv (bv)^{-1}((\al -
\be)\frac{dt}{t} + adt) \\ 
tdt \equiv \frac{bv}{a}\frac{dt}{t} - \frac{\al -\be +1}{a}dt
\end{gather*}
getting finally
\begin{gather*}\nabla_{GM}(\frac{dt}{t}) = \frac{dt}{t}\otimes
((\al - \be)\frac{db}{b} + 
\al \frac{dv}{v}) + dt\otimes (da + a\frac{db}{b} + a\frac{dv}{v}) \\
\nabla_{GM}(dt) = dt\otimes (\be \frac{dv}{v} - (\al -\be +1)\frac{da}{a}) +
\frac{dt}{t}\otimes (bv\frac{da}{a}+d(bv)). 
\end{gather*}
For convenience define
$$\theta = \frac{da}{a}+\frac{db}{b}+\frac{dv}{v}.
$$
Representing an element in our rank two bundle as a column vector
$$\begin{pmatrix}r \\ s \end{pmatrix} = rdt+s\frac{dt}{t}
$$
the matrix for the connection on $\sE$ becomes
$$A := \begin{pmatrix}\be \ta -(\al +1)\frac{da}{a} - \be\frac{db}{b} & a\ta \\
bv\ta & \al\ta - \al \frac{da}{a} -\be \frac{db}{b}
\end{pmatrix}.
$$
The corresponding connection has a regular singular point at
$v=0$ and an irregular one at 
$v=\infty$. Extending $\sE$ to $\sO^2$ on $\P^1_K$, we can take
$\sD = (0)+2(\infty)$ but 
the matrix
$g$ is not invertible at $\infty$. In order to remedy this, make the base
change $z^{-2} = v$, and 
adjoin to
$K$ the element
$\sqrt{ab}$. Notice the base change modifies the Gau\ss-Manin
determinant computation. Let 
us ignore this for a while and continue with the determinant calculation. 

Define $\ga := \frac{\sqrt{ab}}{z}$. Make the change of basis
\begin{multline*} A_{{\rm new}} = \\
\begin{pmatrix}1 & \frac{-z^2}{b}(\be - \ga) \\
0 & \frac{z^2}{2b}\end{pmatrix}
\begin{pmatrix}\be \ta -(\al +1)\frac{da}{a} -
\be\frac{db}{b} & a\ta \\
\frac{b}{z^2}\ta & \al\ta - \al \frac{da}{a} -\be \frac{db}{b}
\end{pmatrix}\times \\
\begin{pmatrix}1 & 2(\be - \ga) \\
0 & \frac{2b}{z^2}\end{pmatrix} +
\begin{pmatrix}1 & \frac{-z^2}{b}(\be - \ga) \\
0 & \frac{z^2}{2b}\end{pmatrix}\begin{pmatrix}0 & -2d\ga \\
0 & d(\frac{2b}{z^2})\end{pmatrix}.
\end{multline*}
This works out to
$$A_{{\rm new}} = \begin{pmatrix}\ga\ta -(\al +1)\frac{da}{a}
-\be \frac{db}{b} & 
(2\al +2\be +1)\ga \ta - 2\be(\al+1)\ta \\
\frac{\ta}{2} & (\al + \be + 2)\ta - \ga\ta -(\al+1)\frac{da}{a} - \be
\frac{db}{b} \end{pmatrix}.
$$
Here of course $\theta= \frac{da}{a} + \frac{db}{b} -
2\frac{dz}{z}$. 

Note
\begin{multline*}\Tr(A_{{\rm new}}) = (\al+\be+2)\ta
-2(\al+1)\frac{da}{a} -2\be\frac{db}{b} \\ 
\equiv (\be
-\al)\frac{da}{a} + (\al
-\be)\frac{db}{b}-2(\al+\be)\frac{dz}{z} \mod d\log(K^\times). 
\end{multline*}
At $z=0$ the polar part of $A_{{\rm new}}$ looks like
$$\begin{pmatrix}\frac{\sqrt{ab}}{z}(\frac{da}{a}+\frac{db}{b}-2\frac{dz}{z}) &
(2\be + 2\al +1)\frac{\sqrt{ab}}{z}(\frac{da}{a}+\frac{db}{b}-2\frac{dz}{z})
+4\be(\al+1)\frac{dz}{z} \\
\frac{-dz}{z} & \frac{-\sqrt{ab}}{z}(\frac{da}{a}+\frac{db}{b}-2\frac{dz}{z}) - 2(\al + \be
+2)\frac{dz}{z}
\end{pmatrix}.
$$
Writing $A_{{\rm pol},0} = g_0\frac{dz}{z^2}+\frac{\eta_0}{z}$, the
matrix for $g_0$ with coefficients in $\C[z]/(z^2)$  is
$$g_0 = \begin{pmatrix}-2\sqrt{ab} & -2(2\al+2\be+1)\sqrt{ab}+4\be(\al+1)z \\
-z & 2\sqrt{ab}-2(\al+\be+2)z\end{pmatrix}.
$$
Also
$$\eta_0 = \begin{pmatrix}\sqrt{ab}(\frac{da}{a}+\frac{db}{b}) &
(2\al+2\be+1)\sqrt{ab}(\frac{da}{a}+\frac{db}{b}) \\
0 & -\sqrt{ab}(\frac{da}{a}+\frac{db}{b})\end{pmatrix}.
$$
With respect to the trivialization $d\log(z^{-1})$ the matrix
$g$ at $z=\infty$ is 
$$g_\infty= \begin{pmatrix}0 & 4\be(\al+1) \\ -1 & -2(\al+\be+2)\end{pmatrix}.
$$
Notice that the matrices for $g$ are invertible both at $0$ and
$\infty$.  
Writing $A_{{\rm pol}, \infty}=
g_\infty\frac{dz}{z}+\eta_\infty$ 
the contribution $\rt dg_\infty g_\infty^{-1} \eta_\infty$ is of
course vanishing, as well as the contribution at $\infty$
obtained by changing the trivialization $\frac{dz}{z}$ to $ {\rm
unit}\cdot \frac{dz}{z}$.
At $0$ we get ($\bar g := g \mod(z)$)
\begin{gather} \label{expr}
\rt (dgg^{-1}\frac{\eta}{z})=\rt(d{\bar g}{\bar g}^{-1}\frac{\eta}{z}) 
=\frac{1}{2}(\frac{da}{a}+\frac{db}{b})\times \\
\Tr
\Big[\begin{pmatrix}0 & 4\be(\al+1) \notag \\
-1 & -2(\al+\be+2) \end{pmatrix} 
\begin{pmatrix}-1 & -(2\be+2\al+1) \notag \\
0 & 1\end{pmatrix}^{-1}
\begin{pmatrix}1 & 2\be+2\al+1 \notag \\
0 & 1\end{pmatrix}\Big]  \notag \\
= \frac{1}{2}(\frac{da}{a}+\frac{db}{b})(2\be+2\al+4) \equiv
(\al+\be)(\frac{da}{a}+\frac{db}{b}) \mod d\log(K^\times).\notag
\end{gather}
Now we compare with the conjectural formula \ref{conj2.7}
$$\det(H^*_{DR}(\sE))^{-1} = c_1(\omega(\sD),s)\cdot \det(\nabla) - 
\rt (dgg^{-1}\frac{\eta}{z}).
$$
Here $s$ can be taken to be the divisor defined by the
trivializing section $\frac{dz}{z^2}-\frac{dz}{z}$ of the sheaf
$\omega(2\cdot 0 +\infty)$, that is $z=1$. Thus one has
$$ \Tr A_{{\rm new}}= (\beta -\alpha)\frac{da}{a} +(\alpha
-\beta)\frac{db}{b}. $$
Further, we have to write
$$s= \frac{dz}{z^2} w,$$
where $w=1-z\in \sO^\times_{X, (20)}$. Since $\Tr \eta_0=0$,
$\rt \frac{dw}{w}\frac{\eta_0}{z}= 0$ as well, thus
the local contribution at $0$ is given by \eqref{expr}.

The conjecture gives (writing $\delta_2 : \P^1 \to \P^1,\
x\mapsto x^2$)
\begin{multline*}\det(H^*_{DR}(\delta_2^*\sE))^{-1} \stackrel{?}{=}
(\be-\al)\frac{da}{a}+(\al-\be)\frac{db}{b}-(\al+\be)(\frac{da}{a}+\frac{db}{b}) \\
 =
-2\al\frac{da}{a}-2\be\frac{db}{b}.
\end{multline*}

Notice we have adjoined $\sqrt{ab}$ to $K$ so we have lost some
$2$-torsion. Bearing in 
mind that $\sE=R^1\pi_{*,DR}$ which introduces a minus sign in
the determinant 
calculations and comparing with our earlier calculation
\eqref{1} above, we find that what 
we need to finish is
\begin{prop} The Gau\ss-Manin determinant for de Rham cohomology
of $\delta_2^*\sE$ is 
twice the corresponding determinant for $\sE$. 
\end{prop}
\begin{proof}Again we use sheaf notation for working with
modules. Recall for the pullback 
we substituted $z^2=w=v^{-1}$. We can write
$\delta_2^*\sE=\sE\oplus z\sE$. We have 
$$\delta_2^*\nabla(ze)=z\nabla(e)+ze\otimes\frac{dv}{-2v},
$$
so with respect to the above decomposition we can write 
$$(\delta_2^*\sE,\delta_2^*\nabla)) = (\sE,\nabla)\oplus(\sE,\nabla
-\frac{1}{2}\frac{dv}{v})
$$
The second term on the right is the connection obtained by tensoring $\sE =
R^1\pi_{*,DR}(\sL_1\boxtimes\sL_2)$ with
$(\sO,-\frac{1}{2}\frac{dv}{v})$. Using the 
projection formula and the invariance of the latter connection,
this is the same as the 
connection on
$R^1\pi_{*,DR}
((\sL_1-\frac{1}{2}\frac{dt}{t})\boxtimes(\sL_2-\frac{1}{2}\frac{du}{u}))$,
i.e. it amounts to replacing $\al,\be$ by
$\al-\frac{1}{2},\be-\frac{1}{2}$. Using (1), 
this changes the Gau{\ss}-Manin 
determinant by $d\log(\sqrt{ab})$ which is trivial. It follows that the
Gau\ss-Manin determinant of 
$\delta_2^*\sE$ is twice that of $\sE$, which is what we want.
\end{proof}
This concludes the proof of Theorem \ref{thm4.3}. 
\end{proof}

\section{Periods}

Let $X/\C$ be a smooth, complete curve. We consider a connection
(relative to $\C$) 
$\nabla:E \to E\otimes\omega_X(\sD)$. Let $\sE$ be the
corresponding local system on 
$X(\C)-\sD$. Notice that we do not assume $\nabla$ has regular
singular points, so 
$\sE$ does not determine $(E,\nabla)$. For example, it can
happen that $\sE$ is a 
trivial local system even though $\nabla$ is highly nontrivial. In
this section, we 
consider the question of associating periods to $\det
H^*_{DR}(X-\sD,E)$. We work 
with  algebraic de Rham cohomology in order to capture the
irregular structure. The 
first remark is that it should be possible using Stokes
structures \cite{M} to write 
down a homological dual group $H_1(X^*,\sE)$ and perfect
pairings ($\sE^\vee$ is the 
dual local system) 
\eq{5.1}{H_1(X^*,\sE^\vee) \times H^1_{DR}(X-\sD,E) \to \C.
}
Here $X^*$ is some modification of the Riemann surface $X$. (The
point is that e.g. 
in the example we give below the de Rham group can be large
while the local system 
$\sE$ is trivial and $X-D=\A^1$.) Let $F\subset \C$ be a
subfield, and assume we are 
given (i) an $F$-structure on $\sE$, i.e. an $F$-local system $\sE_F$ and an
identification $\sE_F\otimes\C \cong \sE$. (ii) A triple
$(X_0,\sD_0,E_0)$ defined 
over $F$ and an identification of the extension to $\C$ of these data with
$(X,\sD,E)$. When e.g. $(E,\nabla)$ satisfies the condition of proposition
\ref{prop2.3}, one has
\ml{5.2}{\det H^*_{DR}(X-\sD,E) \\
\cong \C\otimes_F\det H^*(X_0,E_0)\otimes\det
H^*(X_0,E_0\otimes\omega(\sD_0))^{-1}  
}
so the determinant of de Rham cohomology gets an $F$-structure,
even if $\nabla$ is not necessary itself defined over $F$.
Of course, (i)
determines an $F$-structure on $H_*(X^*,\sE)$. Choosing bases
$\{p_j\},\ \{\eta_k\}$ compatible with the
$F$-structure and taking the determinant of the matrix of periods $\int_{p_j}\eta_k$
\eqref{5.1} yields an invariant
\eq{5.3}{Per(E_0,\nabla,\sE_F)\in \C^\times/F^\times.
}
(More generally, one can consider two subfields $k,F\subset \C$ with a reduction of
$\sE$ to $F$ and a reduction of $E$ to $k$. The resulting determinant lies in
$F^\times\backslash\,\C^\times/k^\times$.) In the case of regular singular points
these determinants have been studied in \cite{ST}. 

Notice that the period invariant depends on the choice of an $F$-structure on the
local system $\sE=\ker(\nabla^{{\rm an}})$. 
When $(E,\nabla)$ are ``motivic'', i.e. come
from the de Rham cohomology of a family of varieties over $X$,
the corresponding 
local system of Betti cohomology gives a natural $\Q$-structure on
$\sE$. By a general 
theorem of Griffiths, the connection $\nabla$ in such a case
necessarily has regular 
singular points. In a non-geometric situation,
or even worse, in the irregular case, there doesn't seem to be
any canonical such 
$\Q$ or $F$-structure. For example, the equation $f'-f=0$ has
solution space 
$\C\cdot e^x$. Is the $\Q$-reduction $\Q\cdot e^x$ more natural
than $\Q\cdot 
e^{x+1}$? Of course, in cases like this where the monodromy is
trivial, the choice of 
$\sE_F$ is determined by choosing an $F$-point $x_0\in
X_0-\sD_0$ and taking $\sE_F 
= E_{x_0}$. 

Even if there is no canonical $F$-structure on $\sE$, one may
still ask for a 
formula analogous to conjecture \ref{conj2.7} for
$Per(E_0,\nabla,\sE_F)$. In this 
final section we discuss the very simplest case
\eq{5.4}{X=\P^1,\  \sD = m\cdot \infty,\  E=\sO_X,\ 
\nabla(1)=df=d(a_{m-1}x^{m-1}+\ldots+a_1x). 
}
Period determinants in this case (and more general confluent
hypergeometric cases) 
were computed by a different argument in \cite{T}. We stress
that our objective here 
is not just to compute the integral, but to exhibit the analogy
with formula \eqref{2.7a}. 
We would like ultimately to find a formula for periods
of higher rank 
irregular connections which bears some relation to
conjectures \ref{conj2.7} and
\ref{conj2.11}. 
We consider the situation \eqref{5.4} with $a_{m-1}\neq 0$. Then
$H^0_{DR}=(0)$ and 
$H^1_{DR}\cong \H^1(\P^1,\sO \to \omega(m\cdot \infty))$ has as basis the
classes of $z^idz,\ 0\le i\le m-3$. One has $\sE^\vee = \C\cdot \exp(f(z))$
(trivial local system) so we take the obvious $\Q$-structure with basis $\exp(f)$. We
want to compute the determinant of the period matrix
\eq{5.5}{\Big(\int_{\sigma_i}\exp(f(z))z^{j-1}dz\Big)_{i,j=1,\dotsc,m-2}
}
for certain chains $\sigma_i$ on some $X^*$. Let $S$ be a union of open sectors
about infinity on $\P^1$ where $Re(f)$ is positive, (i.e. $S$ is a union of
sectors of the form (here $N>>1$ and $\epsilon <<1$ are fixed)
\begin{multline*}S_k :=\{re^{i\theta}\ | \ N< r<\infty,\
\frac{-\text{arg}(a_{m-1})+(2k-\frac{1}{2}-\epsilon)\pi}{m-1}\} <
\theta_k < \\
\frac{-\text{arg}(a_{m-1})+(2k+\frac{1}{2}+\epsilon)\pi}{m-1}\}
\end{multline*} so
$X^*:= \P^1-S\sim
\P^1 -
\{p_1,\ldots,p_{m-1}\}$ where the $p_k$ are distinct points. In
particular, $H_1(X^*) = \Z^{m-2}$.
Define  $\sigma_k := \gamma_k -
\gamma_0$, where 
\eq{5.6}{\gamma_k := \Big\{r\exp(i\theta)\ | 0\le r<\infty;\
\theta=\frac{-\text{arg}(a_{m-1})+(2k+1)\pi}{m-1}\Big\}
}
The $\sigma_k,\ 1\le k\le m-2$ form a basis for $H_1(X^*,\Z)$. 

Write $P_{ij} = \int_{\sigma_i}\exp(f(z))z^{j-1}dz$. 
\begin{lem}\label{lem5.1} We have
\ml{5.7}{\det(P_{ij})_{1\le i,j\le m-2} \\
= \int_{\sigma_1\times\cdots\times
\sigma_{m-2}}\exp(f(z_1)+\ldots+f(z_{m-2}))\prod_{i<j}
(z_j-z_i)dz_1\wedge\ldots\wedge dz_{m-2} }
\end{lem}
\begin{proof}The essential point is
the expansion
\eq{5.8}{\prod_{i<j}(z_j-z_i)=\sum_a (-1)^{{\rm
sgn}(a)}z_1^{a(1)-1}z_2^{a(2)-1}\cdots 
z_{m-2}^{a(m-2)-1}, }
where $a$ runs through permutations of $\{1,\dotsc,m-2\}$. 
\end{proof}

We will evaluate \eqref{5.7} by stationary phase considerations
precisely parallel 
to the techniques described in section 2 and \cite{irreg}. Indeed, the
degree $m-2$ part  $J^{m-2}(\P^1, m \cdot \infty)\subset J(\P^1,
m\cdot \infty)$ of the generalized jacobian is simply the
$\sO_{m\cdot\infty}^\times$-torsor $\omega_{m\cdot\infty}^\times
$ of trivializations 
of
$\omega_{m\cdot\infty} =
\omega(m\cdot\infty)/\omega$ modulo multiplication by a constant
in $\C^\times$. Writing $u=z^{-1}$, we may 
identify this torsor with 
\eq{5.9}{\{b_0\frac{du}{u}+\ldots+b_{m-1}\frac{du}{u^m}\,|\, b_{m-1}\neq 0\}.
}
The quotient of such trivializations up to global isomorphism is
\eq{5.10}{\omega_{m\cdot\infty}^\times/\C^\times = 
\{s_{m-1}\frac{du}{u}+\ldots+s_{1}\frac{du}{u^{m-1}}+\frac{du}{u^{m}}
\}=\{(s_{m-1},\dotsc,s_1)\}. }
Let $B\subset \omega_{m\cdot\infty}^\times/\C^\times$ be defined
by $s_{m-1}=0$. Let
$\Gamma(\P^1,\omega(m\cdot\infty))^\times$ denote the space of sections which
generate $\omega(m\cdot\infty)$ at $\infty$. We have
\eq{5.11}{\A^{m-2}\twoheadrightarrow \text{Sym}^{m-2}(\A^1)
\stackrel{div}{\cong} 
\Gamma(\P^1,\omega(m\cdot\infty))^\times/\C^\times \cong B \subset
\omega_{m\cdot\infty}^\times/\C^\times. 
}
Let $z_1,\dotsc,z_{m-2}$ be as in \eqref{5.7}, and add an extra
variable $z_{m-1}$. 
Take
$s_k(z_1,\dotsc,z_{m-1})$ to be the $k$-th elementary symmetric
function, so e.g. 
$s_{m-1} = z_1z_2\cdots z_{m-1}$. We have a commutative diagram
\eq{5.12}{\begin{CD}\A^{m-2} @>>> B \\ 
@VVz_m-1=0V @VVV \\
\A^{m-1}
@>z\mapsto(s_{m-1}(z),\dotsc,s_1(z))>>\omega_{m\cdot\infty}^\times/\C^\times. 
\end{CD}
}

Notice that
\eq{5.13}{\prod_{i<j}(z_j-z_i)dz_1\wedge\ldots\wedge dz_{m-2} =
ds_1\wedge\ldots\wedge ds_{m-2}. 
}
Let $p_k(z_1,\dotsc,z_{m-2}) = z_1^k+\ldots+z_{m-2}^k$ be the
$k$-th power sum (or $k$-th Newton class). 
Define
\eq{5.14}{F(s_1,\dotsc,s_{m-2}) := f(z_1)+\ldots+f(z_{m-2}) =
a_1p_1+\ldots+a_{m-1}p_{m-1}.  
}
Notice that, although the righthand expression makes sense on all of
$\omega_{m\cdot\infty}^\times/\C^\times$, we think of $F$ as defined only on
$B:s_{m-1}=0$. Let
$\Psi$ be the direct image on
$B$ of the chain
$\sigma_1\times\cdots\times\sigma_{m-2}$ on $\A^{m-2}$. The integral \eqref{5.7}
becomes
\eq{5.15}{\int_\Psi \exp(F(s_1,\dotsc,s_{m-2}))ds_1\wedge\ldots\wedge ds_{m-2}. 
}

\begin{lem}Let $b\in B=\text{Sym}^{m-2}(\A^1)$ correspond to the divisor of zeroes of
$df = f'dz = (a_1+2a_2z+\ldots+(m-1)a_{m-1}z^{m-2})dz$. Then $dF$ vanishes at $b$
and at no other point of $B$. 
\end{lem}
\begin{proof} The differential form 
$$\eta :=
a_1dp_1(z_1,\dotsc,dz_{m-1})+\ldots+a_{m-1}dp_{m-1}(z_1,\dotsc,z_{m-1})
$$
on $\omega_{m\cdot\infty}^\times/\C^\times$ is translation invariant. Indeed, to see
this we may trivialize the torsor and take the point $s_1=\ldots=s_{m-1}=0$ to be
the identity. Introducing a formal variable $T$ with $T^m=0$, the group structure is
then given by $s\oplus s'=:s''$ with 
\ml{}{(1-s_1T+\ldots+(-1)^{m-1}s_{m-1}T^{m-1})(1-s_1'T+\ldots+(-1)^{m-1}
s_{m-1}'T^{m-1})
\\ =(1-s_1''T+\ldots+(-1)^{m-1}s_{m-1}''T^{m-1}).
}
Since $-\log(1-s_1T+\ldots+(-1)^{m-1}s_{m-1}T^{m-1}) =
p_1T+\ldots+p_{m-1}T^{m-1}$ 
it follows that the $p_i$ are additive, whence $\eta$ is
translation invariant. 

Note that $dF=\eta|_B$. Define $\pi:\A^1 \to B\subset 
\omega_{m\cdot\infty}^\times/\C^\times$ by $\pi^*p_k =z^k,\ k\le m-2,\
\pi^*s_{m-1}=0$. Then
$\pi^*\eta = df$. In particular, $\pi^*\eta$ vanishes at the
zeroes of $df=f'dz$. It 
follows that since $b=(f')\in \text{Sym}^{m-2}(\A^1)$, we have
$DF|_b = 0$ as well. 
The proof that $b$ is the unique point where $\eta|_B$ vanishes is given in
\cite{irreg}, lemma 3.10. We shall omit it here. 
\end{proof}

Note that 
$$b=(b_1,\dotsc,b_{m-2}) =
(\frac{a_1}{(m-1)a_{m-1}},\frac{2a_2}{(m-1)a_{m-1}},\dotsc,
\frac{(m-2)a_{m-2}}{(m-1)a_{m-1}})
$$ 
in the $s$-coordinate system on $B$. Set $t_i = s_i-b_i$, and write $F(s) =
F(b)+G(t)$, so $G(t_1,\dotsc,t_{m-2})$ has no constant or linear terms.  
\begin{lem} There exists a non-linear polynomial change of variables of the form
$$t_j' = t_j+B_j(t_1,\dotsc,t_{j-1})
$$
such that $B(0,\dotsc,0)=0$ and $G(t) = Q(t')$ where $Q$ is homogeneous of degree
$2$. 
\end{lem}
\begin{proof} The proof is close to \cite{irreg}, lemma 3.10.
We write (abusively) $p_k(s)$ for the power sum
$z_1^k+\ldots+z_{m-1}^k$, taken as a function of the elementary symmetric
functions $s_1,\dotsc,s_k$. The quadratic monomials $s_is_{m-1-i}$ all occur with
nonzero coefficient in $p_{m-1}$. By construction,
$F(s)=a_1p_1(s)+\ldots+a_{m-1}p_{m-1}(s)$ with $a_{m-1}\neq 0$. If we think of $s_i$
as having weight $i$, $p_k(s)$ is pure of weight $k$, so $s_is_{m-1-i}$ occurs with
nonzero coefficient in $F(s)$. Since the weight $m-1$ is maximal, $t_it_{m-1-i}$
will occur with nonzero coefficient in $G(t)$ as well. Thus, we have
$$G(t) = Q(t_1,\dotsc,t_{m-1})+H(t)
$$
where $Q$ is quadratic and contains $t_it_{m-1-i}$ with nonzero coefficient, and $H$
has no terms of degree $<3$. Further, $H$ has no terms of weight $>m-1$. In
particular, the variable $t_{m-2}$ does not occur in $H$. If we replace $t_{m-2}$ by
$t_{m-2}':=t_{m-2}+A_{m-2}(t_1,\dotsc,t_{m-3})$ for a suitable polynomial $A_{m-2}$,
we can eliminate $t_1$ from $H$ completely,
$G(t)=Q(t_1,\dotsc,t_{m-3},t_{m-2}')+\tilde H(t_2,\dotsc,t_{m-3})$. The weight
and degree conditions on $\tilde H$ are the same as those on $H$, so we conclude
that $\tilde H$ does not involve $t_{m-3}$. Also, 
this change does not affect the
monomials
$t_it_{m-1-i}$ in $Q$ for $i\ge 2$. Thus we may write
$$G = (*)t_1t_{m-2}' + (**)t_2t_{m-3}+\tilde Q(t_3,\dotsc,t_{m-4})+\tilde
H(t_2,\dotsc,t_{m-4})
$$
since $(**) \neq 0$, we may continue in this fashion, writing $t_{m-3}' =
t_{m-3}+A_{m-3}(t_1,\dotsc,t_{m-4})$, etc. 
\end{proof}

The constant term can be written
\eq{5.17}{F(b) = \sum_{\substack{\beta\\
f'(\beta)=0}}f(\beta).
}
The nonlinear change of variables $t\mapsto t'$ has jacobian $1$. Also, the
quadratic form above is necessarily nondegenerate 
(otherwise $F$ would have 
more that one critical point). 
One has
\begin{prop}
 \eq{5.18}{\det (P_{ij})_{1\le i,j,\le m-2}=
\prod_{\substack{\beta\\
f'(\beta)=0}}\exp(f(\beta))\int_\Theta
\exp(Q(t))dt_1\wedge\ldots\wedge dt_{m-2}, 
}
where 
\eq{5.19}{Q = t_1^2+\ldots+t_{m-2}^2
} is the standard, nondegenerate quadric on $B$, and $\Theta$ is
some $n-2$-chain on 
$B$. Moreover, $\int_\Theta$ is determined up
to
$\Q^\times$-multiple on purely geometric grounds.
\end{prop}
Since the shape of the integral is obviously 
coming from the change of coordinates $t \mapsto t'$, we have to
understand the meaning of $\int_\Theta$.

Let $W\subset \P^n\times \P^1$ be the family of quadrics over $\P^1$
defined by
\eq{3b}{UQ(S_1,\dotsc,S_{m-2})-VT^2 = 0.
}
Here $S_1,\dotsc,S_{m-2},T$ are homogeneous coordinates on
$\P^{m-2}$, $U,V$ are 
homogeneous coordinates on $\P^1$, and $Q(S_1,\dotsc,S_{m-2})$ is a
nondegenerate quadric. We have Weil divisors $Y:U=T=0;\ Z:Q=T=0$ in $W$.
Note $Y$ and $Z$ are smooth, and $W_{\text{sing}} = Y\cap Z$. Let $\pi:W'
\to W$ be the blowup of
$W$ along the Weil divisor $Y$. Let $Y'\subset W'$ be the exceptional
divisor. 
\begin{lem}\begin{enumerate}
\item[i.] $W'$ is smooth.
\item[ii.] The strict transform $Z'$ of $Z$ in $W'$ is isomorphic to $Z$
and 
$$Y'\cap Z'\subset Y_{\text{smooth}}'.
$$
\end{enumerate}
\end{lem}
\begin{proof}Let $P'=\text{BL}(Y\subset \P^n\times\P^1)$ be the blowup.
Then $W'$ is the strict transform of $W$ in $P'$. Since $Y\cap Z$ is the
Cartier divisor $U=0$ in $Z$, it follows that the strict transform of $Z$
in $P'$ or $W'$ is isomorphic to $Z$. 

We consider the structure of $W'$ locally around the exceptional divisor.
We may assume some $S_i$ is invertible and write $s_j = S_j/S_i,\
t=T/S_i,\ u=U/V,\ q=Q/S_i^2$. The local defining equation for $W$ is
$uq(s)-t^2=0$. Thinking of $W'$ as $\text{Proj}(\bigoplus_{p\ge 0}
I^p\sO_W)$ with
$I=(u,t)$, we have open sets $\sU_1 : \tilde t\neq 0$ and
$\sU_2: \tilde u\neq 0$. (The tilde indicates we view these as projective
coordinates on the Proj.) We have the following coordinates and equations
for
$W'$ and
$Y'$:
\begin{gather}\label{4c}\sU_1; \qquad u't=u; \qquad W':u'q(s)-t=0; \qquad
Y':t=0.
\\
\sU_2; \qquad t'u=t; \qquad W':q(s)-ut'{}^2=0;\qquad Y':u=0.\notag
\end{gather}
The strict transform $Z'$ of $Z$ lies in the locus $\tilde t=0$ and so
doesn't meet $\sU_1$. Both defining equations for $W'$ are smooth, and
$Y'$ is smooth on $\sU_2$. Finally, $Z'\cap Y':q=t'=u=0$ is also smooth.  
\end{proof}

Write 
$$W^0 := W'-Z';\quad Y^0 = Y'-Y'\cap Z'.
$$
We want to show that the chains over which we integrate
can be understood as chains on the topological pair
$(W^0-U,Y^0)$ for some open $U$
(cf. lemma \ref{lem5.5} below). In
$z$-coordinates, we deal with chains $\gamma_k$, \eqref{5.6},
which are parametrized 
$\alpha_k =  r_ke^{i\theta_k}, 0\le r_k<\infty$ for fixed
$\theta_k$. 
Note that for $r>>1$ the
real part of $f$ on $\gamma_k$ will $\to -\infty$. By abuse of notation,
we write $\gamma_j$ also for the closure of this chain on $\P^1$,
i.e. including the point $r=\infty$. 

Write $F(z) = F(b) +Q(S_1,\dotsc,S_{m-2})$.  It is easy to check by looking at
weights that $|S_k| = O(|r|^k)$ as $|r|\to \infty$. On the other hand, because
the paths are chosen so the real parts of $a_{m-1}\alpha_k^{m-1}$ are all
negative we find there exist positive constants $C,C'$ such that
$C|r|^{m-1}\le |F(z)|=|Q(S(z))|\le C'|r|^{m-1}$. In
homogeneous coordinates
$$(S_k,T), (U,V),
$$
the point associated to a point with coordinates
$z$ on our chain is
\begin{multline}\label{7b}  S_k=S_k(z)=O(|r|^k),\ 1\le k\le m-2;\
T=1;\\  U=1;\  V=Q(S(z))\ge C|r|^{m-1}. 
\end{multline}
With reference to the coordinates in \eqref{4c} we see that
\begin{gather}|u|=|U/V| \le C^{-1}|r|^{1-m},\ |t|=|T/S_i| \ge
C_1|r|^{-i},\\ 
|t'|=|t/u|\ge C_2|r|^{m-1-i}. \notag
\end{gather}
In particular, the limit as $|r|\to \infty$ does not lie on $\sU_2$.
Since, near $\infty$ $Z\subset \sU_2$, we conclude our chains stay away
from $Z$ at infinity. 

We fix $\epsilon <<1$ and $N>>1$ and define a connected, simply connected
domain $D\subset \A^1\subset \P^1$ by
$$D=\{re^{i\theta}\ | \ r> N,\ -\frac{\pi}{2}-\epsilon < \theta <
\frac{\pi}{2}+\epsilon\}
$$
thus, $D$ is an open sector at infinity, and $\exp(z)$ is rapidly
decreasing as $|z|\to \infty$ in the complement 
of $D$. In what follows, let $g:W
\to \P^1$ be the projection.
\begin{lem}\label{lem5.5} The assignment
$$\gamma \mapsto \int_\gamma
\exp(Q(S_1/T,\dotsc,S_{m-2}/T)d(S_1/T)\wedge\ldots\wedge d(S_{m-2}/T)
$$
defines a functional $H_{m-2}(W^0-g^{-1}(D),Y^0;\Q) \to \C$. 
\end{lem} 
\begin{proof}Write $\tau$ for the above integrand. Let $M$ be some
neighborhood of $Z$. Then
$\tau$ is rapidly decreasing on $W^0-g^{-1}(D)$ near $Y^0-M\cap Y^0$,
where the size is defined by some metric on the holomorphic $m-2$ forms on
$W^0$.  Since the  chains are compact, a chain $\gamma$ on $W^0$ will be
supported on $W-M$ for a sufficiently small neighborhood $M$ of $Z$.
Thus, integration defines a functional
$$C_{m-2}(W^0-g^{-1}(D),Y^0) \to \C.
$$
It remains to show $\int_{\partial\Gamma} \tau = 0$ for an $m-1$ chain
$\Gamma$. Let $M$ be an open neighborhood of $Z$ not meeting $\Gamma$.
Let $R$ be an open neighborhood of $Y^0$ in $W-M$. Write $\Gamma =
\Gamma_1+\Gamma_2$ where $\Gamma_1\subset \bar{R}$ and $\Gamma_2\cap R =
\emptyset$. Since $\tau$ is closed, $\int_{\partial\Gamma_2}\tau = 0$. On
the other hand, 
the volume of $\partial\Gamma_1$ can be
taken to be bounded independent of $R$.
It follows that
$\int_{\partial\Gamma} \tau = 0$.
\end{proof}
\begin{lem}$H_{m-2}(W^0-g^{-1}(D),Y^0;\Q)\cong \Q$.
\end{lem}
\begin{proof}Let $p\in D$ be a point. We first show
\eq{9c}{H_{m-2}(W^0-g^{-1}(p),Y^0;\Q)\cong \Q.
}
We calculate $H^*(W^0-g^{-1}(p))$ using the Leray
spectral sequence. For $x\neq 0,\infty$, $g^{-1}(x)$ is a smooth, affine
quadric of dimension $m-3$. So $R^pg_*\Q=(0)$ away from $0,\infty$ for
$p\neq 0, m-3$, and $R^{m-3}g_*\Q|_{\P^1-\{0,\infty\}}$ is a rank $1$
local system.  The monodromy about $0$ and $\infty$ is induced by
$(S,T)\mapsto (S,-T)$ and $(S,T)\mapsto (-S,T)$, respectively. Both
actions give $-1$ on the fibres. It follows that, writing
$j:\P^1-\{0,\infty\}\inj \P^1$, we have
\eq{4b}{j_!R^{m-3}g_*\Q|_{\P^1-\{0,\infty\}}\cong
j_*R^{m-3}g_*\Q|_{\P^1-\{0,\infty\}}\cong
Rj_*R^{m-3}g_*\Q|_{\P^1-\{0,\infty\}}. }

It follows that the natural map
$R^{n-1}g_*\Q\to j_*R^{n-1}g_*\Q|_{\P^1-\{0,\infty\}}$ is surjective and
we get a distinguished triangle in the derived category
\eq{11c}{\sP \to R^{n-1}g_*\Q\to Rj_*R^{n-1}g_*\Q|_{\P^1-\{0,\infty\}}
}
where $\sP$ is a sheaf supported over $0,\infty$. In particular,
\eq{12c}{H^1(\P^1,R^{m-3}g_*\Q) \cong
H^1(\P^1-\{0,\infty\},R^{n-1}g_*\Q)=(0),
}
where the vanishing comes by identifying with group cohomology of $\Z$
acting on $\Q$ with the generator acting by $-1$. An easy Gysin argument
yields
\eq{13c}{H^1(\P^1-\{p\},R^{m-3}g_*\Q) \cong
H^1(\P^1-\{0,\infty,p\},R^{m-3}g_*\Q)=\Q.  
}
It also follows from
\eqref{11c} that $H^2(\P^1,R^{m-3}g_*\Q)=(0)$. The spectral sequence thus
gives
\eq{14c}{ H^{m-2}(W^0,\Q) \cong (R^{m-2}g_*\Q)_{\{0,\infty\}}.
}

To compute these stalks, note the fibre of $W'\to \P^1$ over $0$ is a
singular quadric with singular point $S_1=\ldots=S_{m-2}=0,\ T=1$ away
from $Z$. Thus the fibre of $g: W^0 \to \P^1$ over $0$ is the homogeneous
affine quadric $Q(S)=0$ which is contractible. Further, because $Z$ meets
the fibre of $W'$ smoothly, one has basechange for the non-proper map
$g$, so $(R^{m-2}g_*\Q)_0 = (0)$. At infinity, we have seen again that
$Z$ meets the fibre smoothly, so again one has basechange for $g$. Let
$h:W^0-Y^0 \inj W^0$. It follows that 
\eq{15c}{(R^{m-2}g_*(h_!\Q))_{\{0,\infty\}}
= (0).
}  
Combining \eqref{13c},\eqref{14c},\eqref{15c} yields \eqref{9c}. 

To finish the proof of the lemma, we must show the inclusion
$$(W^0-g^{-1}(D),Y^0)\inj (W^0-g^{-1}(p),Y^0)
$$
is a homotopy equivalence. We can define a homotopy from $D-\{p\}$ to
$\partial D$ by flowing along an outward vector field $v$. E.g. if $p=0$
and one has cartesian coordinates $x,y$, one can take
$v=x\frac{d}{dx}+y\frac{d}{dy}$. Since $W'/\P^1$ is smooth over $D$, one
can lift $v$ to a vector field $w$ on $g'{}^{-1}(D)$. Since $Z'$ meets
the fibres of $g'$ smoothly over some larger $D_1\supset D$, we can
arrange for $w$ to be tangent to $Z'$ along $Z'$. Let $h$ be a smooth
function on $\P^1$ which is positive on $D$ and vanishes on $\P^1-D$.
We view $g'{}^*(h)w$ as a vector field on $W'-g'{}^{-1}(p)$. Flowing
along $g'{}^*(h)w$ lifts the flow along $hv$, carries
$g'{}^{-1}(\bar{D})$ into $g'{}^{-1}(\partial D)$ and stabilizes
$W'-Z'$ over $\bar{D}$. This is the desired homotopy equivalence. 
\end{proof}

\begin{remark}It follows from theorem 2.3.3 in \cite{T} that 
$$\int_\Theta \exp(Q)dt \in
\Big(\frac{2\pi}{(m-1)a_{m-1}}\Big)^\frac{m-2}{2}\cdot\Q. 
$$
\end{remark}

\newpage
\bibliographystyle{plain}
\renewcommand\refname{References}

\end{document}